\documentclass[amstex,12pt,english,amssymb]{article}

\usepackage{mathtext}
\usepackage[cp1251]{inputenc}
\usepackage[T2A]{fontenc}
\usepackage[english]{babel}
\usepackage[dvips]{graphicx}
\usepackage{amsmath}
\usepackage{amssymb}
\usepackage{amsxtra}
\usepackage{latexsym}
\usepackage{ifthen}

\usepackage[dvips]{graphicx}
\usepackage{epstopdf}
\usepackage{epsfig}
\usepackage{epic}
\usepackage{eepic}

\begin{document}

%\markboth{\centerline{D. KOVTONYUK, I. PETKOV AND V. RYAZANOV}}
%{\centerline{On regular homeomorphisms of the class $W^{1,1}_{\rm
%loc}$ on plane}}

\newcounter{lemma}[section]
\newcommand{\lemma}{\par \refstepcounter{lemma}%
{\bf Lemma \arabic{section}.\arabic{lemma}.}}
\renewcommand{\thelemma}{\thesection.\arabic{lemma}}

\newcounter{corollary}[section]
\newcommand{\corollary}{\par \refstepcounter{corollary}%
{\bf Corollary \arabic{section}.\arabic{corollary}.}}
\renewcommand{\thecorollary}{\thesection.\arabic{corollary}}

\newcounter{remark}[section]
\newcommand{\remark}{\par \refstepcounter{remark}%
{\bf Remark \arabic{section}.\arabic{remark}.}}
\renewcommand{\theremark}{\thesection.\arabic{remark}}

\newcounter{theorem}[section]
\newcommand{\theorem}{\par \refstepcounter{theorem}%
{\bf Theorem \arabic{section}.\arabic{theorem}.}}
\renewcommand{\thetheorem}{\thesection.\arabic{theorem}}

\newcounter{proposition}[section]
\newcommand{\proposition}{\par \refstepcounter{proposition}%
{\bf Proposition \arabic{section}.\arabic{proposition}.}}
\renewcommand{\theproposition}{\thesection.\arabic{proposition}}

\renewcommand{\theequation}{\arabic{section}.\arabic{equation}}

\def\kohta #1 #2\par{\par\noindent\rlap{#1)}\hskip30pt
\hangindent30pt #2\par}
\def\A{{{\cal {A}}}}
\def\C{{{\Bbb C}}}
\def\R{{{\Bbb R}}}
\def\Rn{{{\Bbb R}^n}}
\def\lC{{\overline {{\Bbb C}}}}
\def\lRn{{\overline {{\Bbb R}^n}}}
\def\lRm{{\overline {{\Bbb R}^m}}}
\def\lRk{{\overline {{\Bbb R}^k}}}
\def\lBn{{\overline {{\Bbb B}^n}}}
\def\Bn{{{\Bbb B}^n}}
\def\B{{{\Bbb B}}}
\def\dist{{{\rm dist}}}
\let\text=\mbox

\def\Xint#1{\mathchoice
   {\XXint\displaystyle\textstyle{#1}}%
   {\XXint\textstyle\scriptstyle{#1}}%
   {\XXint\scriptstyle\scriptscriptstyle{#1}}%
   {\XXint\scriptscriptstyle\scriptscriptstyle{#1}}%
   \!\int}
\def\XXint#1#2#3{{\setbox0=\hbox{$#1{#2#3}{\int}$}
     \vcenter{\hbox{$#2#3$}}\kern-.5\wd0}}
\def\dashint{\Xint-}

\def\cc{\setcounter{equation}{0}
\setcounter{figure}{0}\setcounter{table}{0}}

\overfullrule=0pt

\title{{\bf On equicontinuity of
homeomorphisms \\ with finite distortion in the plane}}

\author{{\bf T. Lomako, R. Salimov and E. Sevostyanov}\\}
\date{\today \hskip 4mm ({\tt LSS271010.tex})}
\maketitle

\begin{abstract}
It is stated equicontinuity and normality of families
$\frak{F}^{\Phi}$ of the so--called homeomorphisms with finite
distortion on conditions that $K_{f}(z)$ has finite mean
oscillation, singularities of logarithmic type or integral
constraints of the type $\int\Phi\left(K_{f}(z)\right)dx\,dy<\infty$
in a domain $D\subset{\C}.$ It is shown that the found conditions on
the function $\Phi$ are not only sufficient but also necessary for
equicontinuity and normality of such families of mappings.
%It is
%also given applications of these results to families of mappings in
%Sobolev's classes.
\end{abstract}

\large

\section{Introduction}

In the theory of  mappings called quasiconformal in the mean,
conditions of the type
\begin{equation}\label{eq2} \int\limits_{D} \Phi
(Q(z))\ dx\,dy\  <\ \infty\end{equation} are standard for various
characteristics $Q$ of these mappings, see e.g. \cite{Ah},
\cite{Bi}, \cite{Gol}, \cite{GMSV}, \cite{Kr$_1$}--\cite{Ku},
\cite{Per}, \cite{Pes}, \cite{Rya} and \cite{UV}. The study of
classes with the integral conditions (\ref{eq2}) is also actual in
the connection with the recent development of the theory of
degenerate Beltrami equations and the so--called mappings with
finite distortion, see e.g. related references  in the monographs
\cite{IM$_1$} and \cite{MRSY}.\medskip

In the present paper we study the problems of equicontinuity and
normality for wide classes of the homeomorphisms with finite
distortion on conditions that $K_{f}(z)$ has finite mean
oscillation, singularities of logarithmic type or integral
constraints of the type (\ref{eq2}) in a domain $D\subset{\C}.$
\medskip

The concept of the generalized derivative was introduced by
Sobolev in \cite{So}. Given a domain $D$ in the complex plane $\C$
the {\bf Sobolev class} $W^{1,1}(D)$ consists of all functions
$f:D\to\C$ in $L^1(D)$ with first partial generalized derivatives
which are integrable in $D$. A function $f:D\to\C$ belongs to
$W^{1,1}_{\mathrm{loc}}(D)$ if $f\in W^{1,1}(D_*)$ for every open
set $D_*$ with compact closure $\overline{D_*} \subset D$.

Recall that a homeomorphism $f$ between domains $D$ and $D'$ in
$\C$ is called of {\bf finite distortion} if $f\in
W^{1,1}_{\mathrm{loc}}$ and
\begin{equation}\label{eq1.0KR}||f'(z)||^2\leqslant K(z)\cdot
J_{f}(z)\end{equation} with a.e. finite function $K$ where
$||f'(z)||$ denotes the matrix norm of the Jacobian matrix $f'$ of
$f$ at $z\in D$ and $J_{f}(z)=\det f'(z)$, see \cite{IM$_1$}. Later
on, we use the notion $K_{f}(z)$ for the minimal function
$K(z)\geqslant1$ in (\ref{eq1.0KR}). Note that
$||f'(z)||=|f_z|+|f_{\bar{z}}|$ and $J_f(z)=|f_z|^2-|f_{\bar{z}}|^2$
at the points of total differentiability of $f$. Thus,
$K_{f}(z)=\frac{||f'(z)||^2}{|J_{f}(z)|}=\frac{|f_z|+|f_{\bar{z}}|}{|f_z|-|f_{\bar{z}}|}$
if $J_{f}(z)\neq0$, $K_{f}(z)=1$ if $f'(z)=0$, i.e.
$|f_z|=|f_{\bar{z}}|=0$, and $K_{f}(z)=\infty$ at the rest points.

Recall that the {\bf (conformal) modulus} of a family $\Gamma$ of
curves $\gamma$ in ${\C}$ is the quantity
\begin{equation}\label{eq4132}
M(\Gamma)=\inf_{\rho \in \,{\rm adm}\,\Gamma} \int\limits_{{\C}}
\rho ^2 (z)\ \ dx\,dy
\end{equation}
where a Borel function $\rho:{\C}\,\rightarrow [0,\infty]$ is {\bf
admissible} for $\Gamma$, write $\rho \in {\rm adm} \,\Gamma $, if
\begin{equation}\label{eq4133}
\int\limits_{\gamma}\rho \,\,ds\ge 1\ \ \ \ \ \ \ \forall\ \gamma
\in \Gamma\ ,
\end{equation}
where $s$ is a natural parameter of the length on $\gamma$.

One of the equivalent geometric definitions of {\bf
$K-$quasiconformal mappings} $f$ with $K\in [1,\infty)$ given in a
domain $D$ in ${\C}$ is reduced to the inequality
\begin{equation}\label{eq4*}
M(f\Gamma)\le K\,M(\Gamma)
\end{equation}
that holds for an arbitrary family  $\Gamma$  of curves $\gamma$ in
the domain $D$. \medskip

Similarly,  given a domain $D$ in ${\C}$   and a (Lebesgue)
measurable function $Q:D\to[1,\infty]$, a homeomorphism
$f:D\to{\overline{{\C}}}$, ${\overline{{\C}}}={\C}\cup\{\infty\}$,
is called {\bf $Q(z)$ -- homeomorphism} if
\begin{equation} \label{eq19*}
M(f\Gamma )\le \int\limits_D Q(z)\cdot \rho^2 (z)\ \ dx\,dy
\end{equation}
for every family  $\Gamma$  of curves $\gamma$ in $D$ and every
$\rho \in {\rm adm} \,\Gamma $, see e.g. \cite{MRSY}.\medskip

In the case $Q(z)\le K$ a.e., we again come to the inequality
(\ref{eq4*}). In the general case, the latter inequality means that
the conformal modulus of the family $f\Gamma$ is estimated by the
modulus  $M_Q$ of $\Gamma$ with the weight $Q,$
$M(f\Gamma)\le M_{Q}(\Gamma),$
see e.g. \cite{AC$_1$}. The inequality of the type (\ref{eq19*}) was
first stated by O. Lehto and K. Virtanen for quasiconformal mappings
in the plane, see Section V.6.3 in \cite{LV}.

%The relation of the type (\ref{eq19*}) was also stated by K. Bishop,
%V. Gutlyanskii, O. Martio and M. Vuorinen in \cite{BGMV} for
%quasiconformal mappings in space where $Q(z)$ is equal to
%$K_I(z,f).$\medskip

Throughout this paper, $B(z_{0},\,r)=\{ z\in\C: |z_{0}-z|<r\}$,
$S(z_{0},\,r)=\{ z\in\C: |z_{0}-z|=r\}$, $S(r)=S(0,\,r),$
$\mathbb{D}=B(0,\,1)$, $R(r_1,r_2,z_0)$ $=\{ z\,\in\,{\C} :
r_1<|z-z_0|<r_2\} $ and $S^2(x,\,r)=\{ y\in \mathbb{R}^{3}:
|x-y|=r\}$. Let $E,$ $F\subset\overline{{\C}}$ be arbitrary sets.
Denote by $\Gamma(E,F,D)$  a family of all curves
$\gamma:[a,b]\rightarrow \overline{{\C}}$ joining $E$ and $F$ in
$D,$ i.e. $\gamma(a)\in E,\gamma(b) \in F$ and $\gamma(t)\in D$ as
$t \in (a, b).$

The following notion generalizes and localizes the above notion of a
$Q$--ho\-me\-o\-mor\-phism. It is motivated by the ring definition
of Gehring for qua\-si\-con\-for\-mal mappings, see e.g.
\cite{Ge$_3$}, introduced first in the plane, see \cite{RSY$_3$},
and extended  later on to the space case in \cite{RS}, see also
Chapters 7 and 11 in \cite{MRSY}.

Given a domain $D$ in ${\C }$ a (Lebesgue) measurable function
$Q:D\rightarrow\,[0,\infty]$, $z_0\in D,$
a homeomorphism $f:D\rightarrow \overline{{\C}}$ is said to be a
{\bf ring $Q$--homeomorphism at the point $z_0$} if
\begin{equation}\label{eq1}
M\left(f\left(\Gamma\left(S_1,\,S_2,\,R(r_1,r_2,z_0)\right)\right)\right)\
\le \int\limits_{R(r_1,r_2,z_0)} Q(z)\cdot \eta^2(|z-z_0|)\ dx\,dy
\end{equation}
for every ring $R(r_1,r_2,z_0)$ and the circles $S_i=S(z_0, r_i)$,
where $0<r_1<r_2< r_0\,\colon =\,{\rm dist}\, (z_0,\partial D),$ and
every measurable function $\eta : (r_1,r_2)\rightarrow [0,\infty ]$
such that
$$\int\limits_{r_1}^{r_2}\eta(r)\ dr\ \ge\ 1\,.$$
$f$ is called a {\bf ring $Q$--homeomorphism in the domain} $D$ if
$f$ is a ring $Q$--ho\-meo\-mor\-phism at every point $z_0\in D$.
Note that, in particular, homeomorphisms $f:D\rightarrow
\overline{{\C}}$ in the class $W_{loc}^{1,2}$ with $K_f(z)\in
L_{loc}^1(D)$ are ring $Q$--ho\-me\-o\-mor\-phisms with
$Q(z)=K_f(z),$ see e.g. Theorem 4.1 in \cite{MRSY}. A  regular
homeomorphism of the Sobolev class $W^{1,1}_{loc}$ in the plane is a
ring Q-homeomorphism with  $Q(z)$ is equal to the so-called
tangential dilatation, see Theorem 3.1. in \cite{Sa}, cf. Lemma
20.9.1 in \cite{AIM}.
\medskip

The notion of ring $Q$--ho\-me\-o\-mor\-phism can be extended in the
natural way to $\infty$.  More precisely, under $\infty\in D
\subseteq \overline{{\C}}$ a homeomorphism $f:D\rightarrow
\overline{{\C}}$  is called a {\bf ring $Q$--ho\-me\-o\-mor\-phism
at ${\bf \infty}$} if the mapping
$\widetilde{f}=f\left(\frac{z}{\,|z|^2}\right)$  is a ring
$Q^{\,\prime}$--ho\-me\-o\-mor\-phism at the origin with
$Q^{\,\prime}(z)=Q\left(\frac{z}{\,|z|^2}\right).$ In other words, a
mapping $f:{\C}\rightarrow \overline{{\C}}$ is a ring
$Q$--ho\-me\-o\-mor\-phism at $\infty$ iff
$$M\left(f\left(\Gamma\left(S(R_1), S(R_2), R(R_1, R_2, 0)\right)\right)\right)\le
\int\limits_{R(R_1, R_2, 0)} Q(w)\cdot
\eta^2\left(|w|\right)du\,dv$$
%\end{equation}
holds for every ring $R(R_1, R_2, 0)$ in $D$ with
$0<R_1<R_2<\infty,$ $S(R_i)$ and for every measurable function $\eta
: (R_1,R_2)\rightarrow [0,\infty ]$ with
%
%
%\begin{equation}\label{eq11}
%
$\int\limits_{R_1}^{R_2}\eta(r)\ dr\ \ge\ 1\,.$

A continuous mappings $\gamma$ of an open subset $\Delta$ of the
real axis $\mathbb{R}$ or a circle into $D$ is called a {\bf dashed
line}, see e.g. 6.3 in \cite{MRSY}. The notion of the modulus of the
family $\Gamma$ of dashed lines $\gamma$ can be given by analogy,
see (\ref{eq4132}).
%Recall that every open set $\Delta$ in $\mathbb{R}$ consists of a
%countable collection of mutually disjoint intervals. This is the
%motivation for the term.
%
%Given a family $\Gamma$ of dashed lines $\gamma$ in complex plane
%$\C$, a Borel function $\varrho:\C\to[0,\infty]$ is called {\bf
%admissible} for $\Gamma$, write $\varrho\in\mathrm{adm}\,\Gamma$,
%if
%\begin{equation}\label{eq1.2KR}\int\limits_{\gamma}\varrho\,ds\ \geqslant\ 1\end{equation}
%for every $\gamma\in\Gamma$. The {\bf (conformal) modulus} of
%$\Gamma$ is the quantity \begin{equation}\label{eq1.3KR}M(\Gamma)\
%=\
%\inf_{\varrho\in\mathrm{adm}\,\Gamma}\int\limits_{\C}\varrho^2(z)\,dx\,dy\,.\end{equation}
We say
that a property $P$ holds for {\bf a.e.} (almost every)
$\gamma\in\Gamma$ if a subfamily of all lines in $\Gamma$ for which
$P$ fails has the modulus zero, cf. \cite{Fu}.
Later on, we also say that a Lebesgue measurable function
$\varrho:\C\to[0,\infty]$ is {\bf extensively admissible} for
$\Gamma$, write $\varrho\in\mathrm{ext\,adm}\,\Gamma$, if
(\ref{eq4132}) holds for a.e. $\gamma\in\Gamma$, see e.g. 9.2 in
\cite{MRSY}.

Given domains $D$ and $D'$ in $\lC=\C\cup\{\infty\}$,
$z_0\in\overline{D}\setminus\{\infty\}$, and a measurable function
$Q:D\to(0,\infty)$, we say that a homeomorphism $f:D\to D'$ is a
{\bf lower Q-homeomorphism at the point} $z_0$ if
\begin{equation}\label{eq1.4KR}M(f\Sigma_{\varepsilon})\ \geqslant\ \inf\limits_{\varrho\in\mathrm{ext\,adm}\,\Sigma_{\varepsilon}}
\int\limits_{D\cap {R (\varepsilon,\,\varepsilon_0,\,z_0)}}
\frac{\varrho^2(z)}{Q(z)}\,dx\,dy\end{equation} for every ring
$R(\varepsilon,\,\varepsilon_0,\,z_0),\,
\varepsilon\in(0,\varepsilon_0),\, \varepsilon_0\in(0,d_0)\,,$ where
$d_0=\sup\limits_{z\in D}\,|z-z_0|\,,$ and $\Sigma_{\varepsilon}$
denotes the family of all intersections of the circles $S(z_0,r),
r\in(\varepsilon,\varepsilon_0)\,,$ with $D$.

The notion can be extended to the case $z_0=\infty\in\overline{D}$
in the standard way by applying the inversion $T$ with respect to
the unit circle in $\lC$, $T(z)=z/|z|^2$, $T(\infty)=0$,
$T(0)=\infty$. Namely, a homeomorphism $f:D\to D'$ is a {\bf lower
$Q$-homeomorphism at} $\infty\in\overline{D}$ if $F=f\circ T$ is a
lower Q$_*$-homeomorphism with $Q_*=Q\circ T$ at $0$. We also say
that a homeomorphism $f:D\to{\lC}$ is a {\bf lower $Q$-homeomorphism
in} $\partial D$ if $f$ is a lower $Q$-homeomorphism at every point
$z_{0}\in\partial D$. Further we show that every homeomorphism of
finite distortion in the plane is a lower $Q$-homeomorphism with
$Q(z)=K_{f}(z)$ and, thus, the whole theory of the boundary behavior
in \cite{KR$_2$}, see also Chapter 9 in \cite{MRSY} can be applied.

%Note that condition (\ref{eq1.4KR}) is equivalent to the inequality
%\begin{equation}M(f\Sigma_{\varepsilon})\ \geq\
%\int\limits_{\varepsilon}^{\varepsilon_0}\frac{dr}{||\,Q||\,_{1}(r)}\,,
%\label{eq8.1.18}
%\end{equation}
%where
%\begin{equation}||\,Q||\,_{1}(r)\ =\
%\int\limits_{D(z_0,r)}Q^{n-1}(z)\ d\A , \label{eq8.1.19}
%\end{equation}
%$d{\cal {A}}$ corresponds to the area of the circle $D(z_0,r)=D \,
%\bigcap \, S(z_0,r)$. Note that the infimum from the right-hand side
%in (\ref{eq1.4KR}) is attained for the function
%\begin{equation}\varrho_0(z)\ =\
%\frac{Q(z)}{||\,Q||_{1}(|\,z|)}\,. \label{eq8.1.20}
%\end{equation}
%Later we often assume that $Q\equiv0$ outside $D$ and take the
%integrals in (\ref{eq8.1.19}) over the whole circle $S(z_0, r)$.

The following term was introduced in \cite{IR}. Let $D$ be a domain
in the complex plane $\mathbb{C}.$ Recall that a function ${\varphi
}:D\rightarrow \mathbb{R}$ has {\bf finite mean oscillation at a
point} $z_{0}\in {D}$ if
\begin{equation}
\overline{\lim\limits_{\varepsilon \rightarrow 0}}\ \ \ \mathchoice
{{\setbox0=\hbox{$\displaystyle{\textstyle -}{\int}$}
\vcenter{\hbox{$\textstyle -$}}\kern-.5\wd0}}
{{\setbox0=\hbox{$\textstyle{\scriptstyle -}{\int}$}
\vcenter{\hbox{$\scriptstyle -$}}\kern-.5\wd0}}
{{\setbox0=\hbox{$\scriptstyle{\scriptscriptstyle -}{\int}$}
\vcenter{\hbox{$\scriptscriptstyle -$}}\kern-.5\wd0}}
{{\setbox0=\hbox{$\scriptscriptstyle{\scriptscriptstyle -}{\int}$}
\vcenter{\hbox{$\scriptscriptstyle -$}}\kern-.5\wd0}}\!\int_{D(z_{0},%
\varepsilon )}|{\varphi }(z)-\overline{{\varphi }}_{\varepsilon
}(z_{0})|\ dxdy\ <\ \infty ,  \label{eq12.2.4}
\end{equation}%
where
\begin{equation}
\overline{{\varphi }}_{\varepsilon }(z_{0})=\mathchoice
{{\setbox0=\hbox{$\displaystyle{\textstyle -}{\int}$}
\vcenter{\hbox{$\textstyle -$}}\kern-.5\wd0}}
{{\setbox0=\hbox{$\textstyle{\scriptstyle -}{\int}$}
\vcenter{\hbox{$\scriptstyle -$}}\kern-.5\wd0}}
{{\setbox0=\hbox{$\scriptstyle{\scriptscriptstyle -}{\int}$}
\vcenter{\hbox{$\scriptscriptstyle -$}}\kern-.5\wd0}}
{{\setbox0=\hbox{$\scriptscriptstyle{\scriptscriptstyle -}{\int}$}
\vcenter{\hbox{$\scriptscriptstyle -$}}\kern-.5\wd0}}\!\int_{D(z_{0},%
\varepsilon )}{\varphi }(z)\ dxdy\ <\ \infty  \label{eq12.2.5}
\end{equation}%
is the mean value of the function ${\varphi }(z)$ over the disk $%
D(z_{0},\varepsilon )$. We also say that a function ${\varphi
}:D\rightarrow \mathbb{R}$ is of {\bf finite mean oscillation in D},
abbr. ${\varphi }\in $ FMO(D) or simply ${\varphi }\in $ FMO, if
${\varphi }$ has a finite mean oscillation at every point $z_{0}\in
{D}$.

\section{Preliminaries}
\setcounter{equation}{0}

Recall that the {\bf spherical (chordal) metric}
$h(z^{\prime},z^{\prime\prime})$ in $\overline{{{\C}}}$ is equal to
$|\pi(z^{\prime})-\pi(z^{\prime\prime})|$ where $\pi$ is the
stereographic projection of $\overline{{{\C}}}$ on the sphere
$S^2(\frac{1}{2}e_{3},\frac{1}{2})$ in ${{\Bbb R}}^{3},$ i.e., in
the explicit form,
$$h(z^{\prime},\infty)=\frac{1}{\sqrt{1+{|z^{\prime}|}^2}}, \ \
h(z^{\prime},z^{\prime\prime})=\frac{|z^{\prime}-z^{\prime\prime}|}{\sqrt{1+{|z^{\prime}|}^2}
\sqrt{1+{|z^{\prime\prime}|}^2}}\,, \ \  z^{\prime}\ne \infty\ne
z^{\prime\prime}\,.$$
The {\bf spherical diameter of a set} $E$ in $\overline{{\C}}$ is
the quantity $h(E)=\sup\limits_{z^{\prime}, z^{\prime\prime}\in E}
h(z^{\prime}, z^{\prime\prime}).$

A family $\frak{F}$ of continuous mappings from $\C$ into
$\overline{\C}$ is said to be a {\bf normal} if every sequence of
mappings $f_m$ in $\frak{F}$ has a subsequence $f_{m_k}$ converging
to a continuous mapping $f:\C \to \overline{\C}$ uniformly on each
compact set $C\subset \C$. Normality is closely related to the
following notion. A family $\frak{F}$ of mappings $f:\C \rightarrow
\overline{\C}$  is said to be {\bf equicontinuous at a point} $z_0
\in \C$ if for every $\varepsilon > 0$ there is $\delta > 0$ such
that $h \left(f(z),f(z_0)\right)<\varepsilon$ for all $f \in
\frak{F}$ and $z \in \C$ with $|z-z_0|<\delta$. The family
$\frak{F}$ is called {\bf equicontinuous} if $\frak{F}$ is
equicontinuous at every point $z_0 \in \C.$ The following version of
the Arzela -- Ascoli theorem will be useful later on, see e.g.
Section 20.4 in \cite{Va}.

\bigskip
\begin{proposition}\label{pr3**!}{\it\, A family $\frak{F}$ of mappings $f:\C\rightarrow
\overline{\C}$ is normal if and only if $\frak{F}$ is
equicontinuous.}
\end{proposition}

\bigskip

For every non-decreasing function $\Phi:[0,\infty ]\to [0,\infty ]
,$ the {\bf inverse function} $\Phi^{-1}:[0,\infty ]\to [0,\infty ]$
can be well defined by setting
\begin{equation}\label{eq5.5CC} \Phi^{-1}(\tau)\ =\
\inf\limits_{\Phi(t)\ge \tau}\ t\ .
\end{equation} As usual, here $\inf$ is equal to $\infty$ if the set of
$t\in[0,\infty ]$ such that $\Phi(t)\ge \tau$ is empty. Note that
the function $\Phi^{-1}$ is non-decreasing, too.\medskip

\begin{remark}\label{rmk3.333} Immediately by the
definition  it is evident  that
\begin{equation}\label{eq5.5CCC} \Phi^{-1}(\Phi(t))\ \le\ t\ \ \ \ \
\ \ \ \forall\ t\in[ 0,\infty ]
\end{equation} with the equality in (\ref{eq5.5CCC}) except
intervals of constancy of the function $\Phi(t)$.
\end{remark}

\medskip
Since the mapping $t\mapsto t^p$ for every positive $p$ is a
sense--preserving homeomorphism $[0, \infty]$ onto $[0, \infty]$ we
may rewrite Theorem 2.1 from \cite{RSY$_1$} in the following form
which is more convenient for further applications. Here, in
(\ref{eq333Y}) and (\ref{eq333F}), we complete the definition of
integrals by $\infty$ if $\Phi_p(t)=\infty ,$ correspondingly,
$H_p(t)=\infty ,$ for all $t\ge T\in[0,\infty) .$ The integral in
(\ref{eq333F}) is understood as the Lebesgue--Stieltjes integral and
the integrals  in (\ref{eq333Y}) and (\ref{eq333B})--(\ref{eq333A})
as the ordinary Lebesgue in\-te\-grals.

\medskip
\begin{proposition} \label{pr4.1aB}{\it\, Let $\Phi:[0,\infty ]\to [0,\infty ]$ be a
non-decreasing function. Set \begin{equation}\label{eq333E} H_p(t)\
=\ \log \Phi_p(t)\ , \qquad \Phi_p(t)=\Phi\left(t^p\right)\,,\quad
p\in (0,\infty)\,.\end{equation}

Then the equality \begin{equation}\label{eq333Y}
\int\limits_{\delta}^{\infty} H^{\,\prime}_p(t)\ \frac{dt}{t}\ =\
\infty
\end{equation} implies the equality \begin{equation}\label{eq333F}
\int\limits_{\delta}^{\infty} \frac{dH_p(t)}{t}\ =\ \infty
\end{equation} and (\ref{eq333F}) is equivalent to
\begin{equation}\label{eq333B} \int\limits_{\delta}^{\infty}H_p(t)\
\frac{dt}{t^2}\ =\ \infty
\end{equation} for some $\delta>0,$ and (\ref{eq333B}) is equivalent to
every of the equalities: \begin{equation}\label{eq333C}
\int\limits_{0}^{\Delta}H_p\left(\frac{1}{t}\right)\ {dt}\ =\ \infty
\end{equation} for some $\Delta>0,$ \begin{equation}\label{eq333D}
\int\limits_{\delta_*}^{\infty} \frac{d\eta}{H_p^{-1}(\eta)}\ =\
\infty
\end{equation} for some $\delta_*>H(+0),$ \begin{equation}\label{eq333A}
\int\limits_{\delta_*}^{\infty}\ \frac{d\tau}{\tau \Phi_p^{-1}(\tau
)}\ =\ \infty \end{equation} for some $\delta_*>\Phi(+0).$
\medskip

Moreover, (\ref{eq333Y}) is equivalent  to (\ref{eq333F}) and hence
(\ref{eq333Y})--(\ref{eq333A})
 are equivalent each to other  if $\Phi$ is in addition absolutely continuous.
In particular, all the conditions (\ref{eq333Y})--(\ref{eq333A}) are
equivalent if $\Phi$ is convex and non--decreasing.}
\end{proposition}

\medskip
It is easy to see that conditions (\ref{eq333Y})--(\ref{eq333A})
become weaker as $p$ increases, see e.g. (\ref{eq333B}). It is
necessary to give one more explanation. From the right hand sides in
the conditions (\ref{eq333Y})--(\ref{eq333A}) we have in mind
$+\infty$. If $\Phi_p(t)=0$ for $t\in[0,t_*]$, then $H_p(t)=-\infty$
for $t\in[0,t_*]$ and we complete the definition $H_p'(t)=0$ for
$t\in[0,t_*]$. Note, the conditions (\ref{eq333F}) and
(\ref{eq333B}) exclude that $t_*$ belongs to the interval of
integrability because in the contrary case the left hand sides in
(\ref{eq333F}) and (\ref{eq333B}) are either equal to $-\infty$ or
indeterminate. Hence we may assume in (\ref{eq333Y})--(\ref{eq333C})
that $\delta>t_0$, correspondingly, $\Delta<1/t_0$ where $t_0\colon
=\sup\limits_{\Phi_p(t)=0}t$, $t_0=0$ if $\Phi_p(0)>0$.

\cc
\section{The main results}

\begin{proposition}
{}\label{prKPR3.1} Let $f:D\to\C$ be a homeomorphism with finite
distortion. Then $f$ is a lower $Q$-homeomorphism at each point
$z_0\in\overline{D}$ with $Q(z)=K_{f}(z)$, see Theorem 3.1. in
\cite{KPR}.
\end{proposition}

\begin{proposition}{}\label{pr8.4.8} Let $D$ and $D'$ be domains in $\C$, let
$z_0\in\overline{D}\setminus\{\infty\}$, and let $Q:D\to(0,\infty)$
be a measurable function. A homeomorphism $f:D\to D'$ is a lower
$Q$-homeomorphism at $z_0$ if and only if
\begin{equation}
M(f\Sigma_{\varepsilon})\ \geq\
\int\limits_{\varepsilon}^{\varepsilon_0}
\frac{dr}{||\,Q||\,_{1}(r)}\quad\forall\
\varepsilon\in(0,\varepsilon_0)\,,\quad\varepsilon_0\in(0,d_0),
\label{eq8.4.9}
\end{equation}
where
\begin{equation}
d_0\ =\ \sup\limits_{z\in D}\, |z-z_0|, \label{eq8.4.10}
\end{equation}
$\Sigma_{\varepsilon}$ denotes the family of all the intersections
of the circles $S(z_0,\,r)$, $r\in(\varepsilon,\varepsilon_0)$, with
$D$, and
\begin{equation}
||\,Q||\,_{1}(r)=\int\limits_{D(z_0,r)}Q(z)\ ds \label{eq8.4.11}
\end{equation} is the
$L_{1}$-norm of $Q$ over $D(z_0,r)=\{z\in D:|\,z-z_0|=r\}=D\cap
S(z_0,r)$. The infimum of the expression from the right-hand side in
(\ref{eq1.4KR}) is attained only for the function
$$\varrho_0(z)\ =\ \frac{Q(z)}{||\,Q||_{1}(|\,z|)}\,,$$
see Theorem 2.1 in \cite{KR$_2$}.
\end{proposition}

\begin{proposition}
{}\label{pr6.4.10} Let $D$ be a domain in $\C$ and
$Q:D\rightarrow \lbrack 0,\infty ]$ a measurable function. A homeomorphism $%
f:D\rightarrow {\C}$ is a ring $Q$-homeomorphism at a point $%
z_{0}$ if and only if, for every $0<r_{1}<r_{2}<d_{0}=\mathrm{dist}%
\,(z_{0},\partial D),$
\begin{equation}
M({\Delta }(fS_{1},fS_{2},fD))\ \leq \ \frac{2\pi}{I},
\label{eq6.3.111}
\end{equation}%
where ${\omega }$ is the area of the unit circle in $\C,$
$q_{z_{0}}(r)$ is the mean value of $Q(z)$ over the circle $|z-z_{0}|=r,$ $%
S_{j}=S(z_{0},\,r_{j}),$ $j=1,2,$ and
$$
I\ =\ I(r_{1},r_{2})\ =\ \int\limits_{r_{1}}^{r_{2}}\
\frac{dr}{rq_{z_{0}}(r)}
$$

Moreover, the infimum from the right-hand side in (\ref{eq1}) holds
for the function
\begin{equation}
\eta _{0}(r)=\frac{1}{Irq_{z_{0}}(r)}\ , \label{eq6.3.116}
\end{equation}
see Theorem 3.15 in \cite{RS}.
\end{proposition}

The above results now yield the following.

\begin{lemma}
{}\label{lem6.3.2} Let $D$ and $D'$ be domains in $\C$, let
$z_0\in\overline{D}\setminus\{\infty\}$, and let $Q:D\to(0,\infty)$
be a measurable function. A homeomorphism $f:D\to D'$ is a lower
$Q$-homeomorphism at $z_0$. Then $f$ is a ring $Q$--homeomorphism at
$z_0$.
\end{lemma}

{\it Proof of Lemma \ref{lem6.3.2}.} Denote by
$\Sigma_{\varepsilon}$ the family of all circles $S(z_0,\,r)$,
$r\in(\varepsilon,\varepsilon_0)$, $\varepsilon_0\in(0,d_0)\,.$ By
Theorem 3.13 in \cite{Zi}, we have
\begin{equation}
M\left(\Delta\left(fS_{\varepsilon},\,fS_{\varepsilon_{0}},\,f(D)\right)\right)\leq
\frac{1}{M\left(f
\Sigma_{\varepsilon}\right)}\leq\frac{2\pi}{\int\limits_{\varepsilon}^{\varepsilon_{0}}\frac{dr}{rq_{z_{0}}(r)}}
\label{eq6.3.1162}
\end{equation}
because
$f\Sigma_{\varepsilon}\subset\Sigma\left(fS_{\varepsilon},\,fS_{\varepsilon_{0}}\right)$,
where $\Sigma\left(fS_{\varepsilon},\,fS_{\varepsilon_{0}}\right)$
consists of all closed curves in $f(D)$ that separate
$fS_{\varepsilon}$ and $fS_{\varepsilon_{0}}$.

Proposition \ref{prKPR3.1} and Lemma \ref{lem6.3.2} imply the
following result.

\begin{theorem} \label{th3.1989898}
Let $f:D\to\C$ be a homeomorphism with finite distortion. Then $f$
is a ring $Q$-homeomorphism at each point $z_0\in\overline{D}$ with
$Q(z)=K_{f}(z)$.
\end{theorem}

\section{Estimates of Distortion}

The results of the following section can be obtained on the base of theorem 3.1 and
the correspondent theorems of work \cite{RS}.

\begin{lemma}
{}\label{lem6.3.23} Let $D$ be a domain in ${\C},$
let $D^{\prime }$ be a domain in $\overline{\C}$ with $h(%
\overline{\C}\setminus D^{\prime })\geq {\Delta }>0,$ and let $%
f:D\rightarrow D^{\prime }$ be a homeomorphism with finite
distortion at a point $%
z_{0}\in D.$ If, for $0<\varepsilon
_{0}<\mathrm{dist}(z_{0},\partial D),$
\begin{equation}
\int\limits_{\varepsilon <|z-z_{0}|\ <\ \varepsilon _{0}}K_{f}(z)\cdot \psi _{{%
\varepsilon }}^{2}(|z-z_{0}|)\ dx\,dy\ \leq \ c\cdot
I^{p}(\varepsilon )\ ,\ \ \ \ \ \ \varepsilon \in (0,\varepsilon
_{0}), \label{eq6.3.24}
\end{equation}%
where $p\leq 2$ and $\psi _{{\varepsilon }}(t)$ is nonnegative function on
$(0,\infty )$ such that
\begin{equation}
0\ <\ I(\varepsilon )\ =\ \int\limits_{\varepsilon }^{\varepsilon
_{0}}\psi _{{\varepsilon }}(t)\ dt<\infty ,\ \ \ \ \ \ \varepsilon
\in (0,\varepsilon _{0}),  \label{eq6.3.25}
\end{equation}%
then
\begin{equation}
h(f(z),f(z_{0}))\ \leq \ \frac{32}{{\Delta }}\ \exp
\left\{-\left(\frac{2\pi}{c}\right)I^{2-p}(|z-z_{0}|)\right\}
\label{eq6.3.26}
\end{equation}%
for all $z\in B(z_{0},{{\varepsilon }_{0}}).$

\end{lemma}

\begin{corollary}
\label{cor6.3.28} Under the conditions of Lemma \ref{lem6.3.23} and
for $p = 1$,
\begin{equation}
h(f(z),f(z_{0}))\ \leq \ \frac{32}{{\Delta }}\ \exp
\left\{-\left(\frac{2\pi}{c}\right)I(|z-z_{0}|)\right\}.
\label{eq6.3.29}
\end{equation}
\end{corollary}

\begin{theorem}
{}\label{th6.4.1} Let $D$ be a domain in ${\C},$ let $%
D^{\prime }$ be a domain in $\overline{\C}$ with $h(\overline{\C}\setminus D^{\prime })\geq {\Delta }>0,$ and let $%
f:D\rightarrow D^{\prime }$ be a homeomorphism with finite
distortion at a point $%
z_{0}\in D.$ Then
\begin{equation}
h(f(z),f(z_{0}))\ \leq \ \frac{32}{{\Delta }}\ \exp \left\{
-\int\limits_{|z-z_{0}|}^{\varepsilon
(z_{0})}\frac{dr}{rq_{z_{0}}(r)}\right\}  \label{eq6.4.2}
\end{equation}%
for $z\in B(z_{0},\varepsilon (z_{0})),$ where $\varepsilon
(z_{0})<\dist(z_{0},\partial D)$ and $q_{z_{0}}(r)$ is the mean
integral value of $K_{f}(z)$ over the circle $|z-z_{0}|=r.$
\end{theorem}

%\begin{remark}
%\label{rem6.4.3} Of course, the mean value $q_{z_{0}}(r)$ of
%$K_{f}(z)$ over some circles $|z-z_{0}|=r$ can be infinite. However,
%$q_{z_{0}}(r)$ is measurable in the parameter $r$ because $K_{f}(z)$
%is measurable in $z,$ say by the Fubini theorem. Moreover, at every
%point $z\neq z_{0},$
%\begin{equation}
%\int\limits_{|z-z_{0}|}^{\varepsilon
%(z_{0})}\frac{dr}{rq_{z_{0}}(r)}\ <\ \infty  \label{eq6.4.4}
%\end{equation}%
%for any a homeomorphism with finite distortion because in the
%contrary case we would have from (\ref{eq6.4.2}) that
%$f(z)=f(z_{0}).$ The integral in (\ref{eq6.4.4})
%can be $0$ if $q_{z_{0}}(r)=\infty $ a.e., but then inequality (\ref{eq6.4.2}%
%) is obvious because ${\alpha }=32$ and ${\Delta }$ as well as $%
%h(f(z),f(z_{0}))$ are less than or equal to $1.$ \medskip
%
%\end{remark}

\begin{corollary}
\label{cor6.4.5} If
\begin{equation}
q_{z_{0}}(r)\leq { \log {\frac{1}{r}} } \label{eq6.4.6}
\end{equation}%
for $r<{\varepsilon }(z_{0})<\mathrm{dist}(z_{0},\partial D),$ then
\begin{equation}
h(f(z),f(z_{0}))\leq \,\frac{32}{{\Delta }}\,\frac{\log \frac{1}{%
\varepsilon _{0}}}{\log \frac{1}{|z-z_{0}|}}  \label{eq6.4.7}
\end{equation}%
for all $z\in B(z_{0},\varepsilon (z_{0}))$.
\end{corollary}

\begin{corollary}
\label{cor6.4.8} If
\begin{equation}
K_{f}(z)\leq { \log {\frac{1}{|z-z_{0}|}} },\,\,\ \text{\ \ }%
\ z\in B(z_{0},\varepsilon (z_{0})),  \label{eq6.4.9}
\end{equation}%
then (\ref{eq6.4.7}) holds in the ball $B(z_{0},\varepsilon
(z_{0}))$.
\end{corollary}

\begin{remark}
\label{rmk6.4.10} If, instead of (\ref{eq6.4.6}) and
(\ref{eq6.4.9}), we have the conditions
\begin{equation}
q_{z_{0}}(r)\leq c\cdot {\log {\frac{1}{r}} } \label{eq6.4.11}
\end{equation}%
and, correspondingly,
\begin{equation}
K_{f}(z)\leq c\cdot { \log {\frac{1}{|z-z_{0}|}} }, \label{eq6.4.12}
\end{equation}%
then
\begin{equation}
h(f(z),f(z_{0}))\ \leq \ {\ \frac{32}{{\Delta }}\ \left[ \frac{\log
\frac{1}{\varepsilon _{0}}}{\log \frac{1}{|z-z_{0}|}}\right] }%
^{1/c}.  \label{eq6.4.13}
\end{equation}
\end{remark}

Choosing in Lemma \ref{lem6.3.23} $\psi (t)=1/t$ and $p=1,$ we also
have the following conclusion.

\begin{corollary}
\label{4.23} Let $f:{\mathbb{D}}\rightarrow {\mathbb{D}}$ be a
homeomorphism with finite distortion such that $f(0)=0$ and
\begin{equation}
\int\limits_{\varepsilon <|z|<1}K_{f}(z)\ \ \frac{dx\,dy}{{|z|}^{2}}\ \leq \ c%
\text{ }\log {\frac{1}{\varepsilon }},\text{ \ }\ \ \ \varepsilon
\in (0,1). \label{eq6.4.24}
\end{equation}%
Then
\begin{equation}
|f(z)|\ \leq \ 64\cdot {|z|}^{\frac{2\pi}{c}}. \label{eq6.4.25}
\end{equation}%

\end{corollary}

\begin{theorem}
{}\label{th6.5.11} Let $D$ be a domain in ${\C},$  let $D^{\prime }$
be a domain in $\overline{\C}$ with $h(\overline{\C}\setminus
D^{\prime })\geq {\Delta }>0,$ and let $f:D\rightarrow D^{\prime }$
be a homeomorphism with finite distortion at a point $z_{0}\in D.$
If $K_{f}(z)$ has finite mean oscillation at the point $z_{0}\in D$,
then
\begin{equation}
h(f(z),f(z_{0}))\leq \frac{32}{{\Delta }}{\left\{ {\frac{\log
\,\frac{1}{\varepsilon _{0}}}{\log \,{\frac{1}{|z-z_{0}|}}}}\right\}
}^{\beta _{0}}  \label{eq6.5.12}
\end{equation}%
for some $\varepsilon _{0}<\mathrm{dist}(z_{0},\partial D)$ and
every $z\in B(z_{0},\varepsilon _{0}),$ where $\beta _{0}>0$ depends
only on the function $K_{f}$.
\end{theorem}

\section{On Normal Families of homeomorphisms with finite
distortion}

The results stated bellow can be proved by theorem 3.1 and the correspondent criteria of normality from the paper \cite{RS}.

Given a domain $D$ in ${\C}$, let $\mathfrak{F}_{K_{f},\Delta }(D)$
be the class of all homeomorphisms $f$ with finite
distortion $K_{f}$ in $D$ with $h(%
\overline{\C}\setminus f(D))\geq {\Delta }>0.$

\begin{theorem}
{}\label{th6.6.1} If $K_{f}\in \mathrm{FMO}$, then
$\mathfrak{F}_{K_{f},\Delta }(D)$ is a normal family.
\end{theorem}

\begin{corollary}
\label{cor6.6.2} The class $\mathfrak{F}_{K_{f},\Delta }(D)$ is
normal if
\begin{equation}
\overline{\lim\limits_{\varepsilon \rightarrow 0}}\ \ \mathchoice
{{\setbox0=\hbox{$\displaystyle{\textstyle -}{\int}$}
\vcenter{\hbox{$\textstyle -$}}\kern-.5\wd0}}
{{\setbox0=\hbox{$\textstyle{\scriptstyle -}{\int}$}
\vcenter{\hbox{$\scriptstyle -$}}\kern-.5\wd0}}
{{\setbox0=\hbox{$\scriptstyle{\scriptscriptstyle -}{\int}$}
\vcenter{\hbox{$\scriptscriptstyle -$}}\kern-.5\wd0}}
{{\setbox0=\hbox{$\scriptscriptstyle{\scriptscriptstyle -}{\int}$}
\vcenter{\hbox{$\scriptscriptstyle -$}}\kern-.5\wd0}}\!\int_{B(z_{0},%
\varepsilon )}K_{f}(z)\ \ dx\,dy\ <\ \infty \ \ \ \ \ \ \forall \
z_{0}\in D. \label{eq6.6.3}
\end{equation}
\end{corollary}

\begin{corollary}
\label{cor6.6.4} The class $\mathfrak{F}_{K_{f},\Delta}(D)$ is normal if every $%
z_0\in D$ is a Lebesgue point of $K_{f}(z)$.
\end{corollary}

\begin{theorem}
{}\label{th6.6.5} Let $\Delta >0$ and let $Q:D\rightarrow \lbrack
0,\infty ]$ be a measurable function such that
\begin{equation}
\int\limits_{0}^{{\varepsilon
}(z_{0})}\frac{dr}{rq_{z_{0}}(r)}=\infty  \label{eq6.6.6}
\end{equation}%
holds at every point $z_{0}\in D$, where ${\varepsilon }(z_{0})=\mathrm{dist}%
(z_{0},\partial D)$ and $q_{z_{0}}(r)$ denotes the mean integral value of $%
K_{f}(z)$ over the circle $|z-z_{0}|=r$. Then
$\mathfrak{F}_{K_{f},\Delta }$ forms a normal family.
\end{theorem}

\begin{corollary}
\label{cor6.6.7} The class $\mathfrak{F}_{K_{f},\Delta }(D)$ is
normal if $K_{f}(z)$ has singularities of the logarithmic type of
order not greater than 1 at every point $z\in D$.
\end{corollary}

%\begin{remark} ? Note that all the above results hold for
%homeo\-mor\-phisms $f$ of the
%Sobolev class $W_{\mathrm{loc}}^{1,2}$ with $f^{-1}\in W_{\mathrm{loc}%
%}^{1,2} $ under the condition that
%\begin{equation}
%K_{I}(z,f)\leq K_{f}(z)\,\,\,\,\,\,\mathrm{a.e.,}  \label{eq6.6.15}
%\end{equation}%
%where $K_{I}$ is the inner dilatation of the mapping $f;$ see, e.g.,
%Theorem 6.1 in \cite{MRSY}. In particular, this is valid for homeomorphisms $f\in W_{%
%\mathrm{loc}}^{1,2}$ with $K_{I}\in L_{\mathrm{loc}}^{1};$ see
%Corollary 6.4 in \cite{MRSY}.
%\end{remark}
%
%\textbf{\textsc{Postscript.}} A family $\mathfrak{F}$ of
%homeomorphisms with finite
%distortion $f:D\rightarrow {\overline{\C}}$ in a domain $%
%D\subset {\C}$  is normal under every condition given above for
%$K_{f}$ if there is a number ${\Delta }>0$ such that one of the
%following conditions holds: \medskip
%
%\noindent \noindent (1) Every mapping $f\in \mathfrak{F}$ omits two
%values whose spherical distance is greater than ${\Delta }.$
%\medskip
%
%\noindent (2) Every mapping $f\in \mathfrak{F}$ omits one value $w_{0}$ and $%
%h\left( w(z_{i}),w_{0}\right) >{\Delta },$ $i=1,2,$ at two fixed points $%
%z_{1}$ and $z_{2}\in D.$ \medskip
%
%\noindent (3) At three fixed points $z_{1},z_{2},$ and $z_{3}\in D,$ $%
%h(w(z_{i}),w(z_{j}))>{\Delta },$ $i\neq j,$ $i,j=1,2,3.$ \medskip
%
%In particular, $\mathfrak{F}$ is normal if all mappings
%$f\in\mathfrak{F}$ omit two fixed values in $\overline{{\C}}.$

\medskip
\section{On some integral conditions}
\setcounter{equation}{0}

The following results can be found in \cite{RS_{2}}.

Recall that a function  $\Phi :[0,\infty ]\to [0,\infty ]$ is called
{\bf convex} if
$$
\Phi (\lambda t_1 + (1-\lambda) t_2)\ \le\ \lambda\ \Phi (t_1)\ +\
(1-\lambda)\ \Phi (t_2)$$ for all $t_1$ and $t_2\in[0,\infty ] $ and
$\lambda\in [0,1]$.\medskip

In what follows, ${\R}(\varepsilon),$ $\varepsilon\in (0,1)$ denotes
the ring in the space ${\C}$,
\begin{equation}\label{eq5.5Cf} {\R}(\varepsilon)\ =R\,(\varepsilon,\,1,\,0).\end{equation}

\noindent The following statement is a generalization and
strengthening of Lemma 3.1 from \cite{RSY$_1$}.

\begin{lemma} \label{lem5.5C} {\it\, Let $Q:{\Bbb D}\to [0,\infty ]$ be a measurable
function and let $\Phi:[0,\infty ]\to (0,\infty ]$ be a
non-decreasing convex function. Suppose that the mean value
$M(\varepsilon)$ of the function $\Phi\circ Q$ over the ring
${\R}(\varepsilon),$ $\varepsilon\in (0, 1),$ is finite. Then
\begin{equation}\label{eq3.222} \int\limits_{\varepsilon}^{1}\
\frac{dr}{rq^{\frac{1}{p}}(r)}\ \ge\ \frac{1}{2}\
\int\limits_{eM(\varepsilon)}^{\frac{M(\varepsilon)}{\varepsilon^2}}\
\frac{d\tau}{\tau \left[\Phi^{-1}(\tau
)\right]^{\frac{1}{p}}}\qquad\qquad\forall\quad p\in (0, \infty)
\end{equation} where $q(r)$ is the average of the function $Q(z)$
over the circle $|z|=r.$ }\end{lemma}

\bigskip

\begin{remark}\label{rmk3.333A} Note that (\ref{eq3.222}) is
equivalent  for each $p\in (0, \infty)$ to the inequality
\begin{equation}\label{eq3.1!}
\int\limits_\varepsilon^1\frac{dr}{rq^{\frac{1}{p}}(r)}\ \ge\
\frac{1}{2}\int\limits_{eM(\varepsilon)}^{\frac{M(\varepsilon)}{\varepsilon^2}}\frac{d\tau}{\tau\Phi_p^{\,-1}(\tau)}\
,\qquad \Phi_p(t)\ \colon =\ \Phi(t^p)\ .
\end{equation} Note also that $M(\varepsilon)$ converges as $\varepsilon\to
0$ to the average of $\Phi\circ Q$ over the unit disk ${\Bbb B}$.
\end{remark}

\medskip
\begin{corollary} \label{cor3.1}{\,\it
Let $\Phi:[0,\infty ]\rightarrow (0,\infty ]$ be a non-decreasing
convex function, $Q:{\Bbb B}\rightarrow [0,\infty ]$ a measurable
function, $Q_*(z)=1$ if $Q(z)<1$ and $Q_*(z)=Q(z)$ if $Q(z)\ge 1$.
Suppose that the mean $M_*(\varepsilon)$ of the function $\Phi\circ
Q$ over the ring ${\R}(\varepsilon),$ $\varepsilon\in (0, 1),$ is
finite. Then
\begin{equation}\label{eq3.1}
\int\limits_{\varepsilon}^{1}\ \frac{dr}{rq^{\frac{\lambda}{p}}(r)}\
\ge\ \frac{1}{2}\
\int\limits_{eM_*(\varepsilon)}^{\frac{M_*(\varepsilon)}{\varepsilon^2}}\
\frac{d\tau}{\tau \left[\Phi^{-1}(\tau )\right]^{\frac{1}{p}}}\
\qquad \  \ \forall\ \lambda\ \in\ (0,1), \qquad p\in (0, \infty)
\end{equation}
where $q(r)$ is the average of the function $Q(z)$ over the circle
$|z|=r.$ }
\end{corollary}

\medskip
\medskip
Indeed, let $q_*(r)$ be the average of the function $Q_*(z)$ over
the circle $|z|=r$. Then $q(r)\le q_*(r)$ and, moreover, $q_*(r)\ge
1$ for all $r\in (0,1)$. Thus, $q^{\frac{\lambda}{p}}(r)\le
q_*^\frac{\lambda}{p}(r)\le q_*^\frac{1}{p}(r)$  for all $\lambda\in
(0,1)$ and hence by Lemma \ref{lem5.5C} applied to $Q_*(z)$ we
obtain (\ref{eq3.1}).

\medskip
\begin{theorem} \label{th5.555}{\it\, Let $Q:{\Bbb D}\to [0,\infty ]$ be a measurable
function such that \begin{equation}\label{eq5.555}
\int\limits_{{\Bbb B}} \Phi (Q(z))\ dx\,dy\  <\ \infty\end{equation}
where $\Phi:[0,\infty ]\to [0,\infty ]$ is a non-decreasing convex
function such that
\begin{equation}\label{eq3.333a} \int\limits_{\delta_0}^{\infty}\ \frac{d\tau}{\tau
\left[\Phi^{-1}(\tau )\right]^{\frac{1}{p}}}\ =\ \infty\,,\qquad
p\in (0, \infty)\,,
\end{equation} for some $\delta_0\
>\ \tau_0\ \colon =\ \Phi(0).$ Then \begin{equation}\label{eq3.333A}
\int\limits_{0}^{1}\ \frac{dr}{rq^{\frac{1}{p}}(r)}\ =\ \infty
\end{equation} where $q(r)$ is the average of the function $Q(z)$
over the circle $|z|=r$.}
\end{theorem}

\begin{remark}\label{rmk4.7www} Since $\left[\Phi^{\,-1}(\tau)\right]^{\frac{1}{p}}=
\Phi_p^{\,-1}(\tau)$ where $\Phi_p(t)=\Phi(t^p),$  (\ref{eq3.333a})
implies that
\begin{equation}\label{eq3.a333} \int\limits_{\delta}^{\infty}\ \frac{d\tau}{\tau
\Phi^{-1}_p(\tau )}\ =\ \infty\ \ \ \ \ \ \ \ \ \ \forall\ \delta\
\in\ [0,\infty)   \end{equation} but (\ref{eq3.a333}) for some
$\delta\in[0,\infty)$, generally speaking, does not imply
(\ref{eq3.333a}). Indeed, for $\delta\in [0,\delta_0),$
(\ref{eq3.333a}) evidently implies (\ref{eq3.a333}) and, for
$\delta\in(\delta_0,\infty)$, we have that
\begin{equation}\label{eq3.e333} 0\ \le\ \int\limits_{\delta_0}^{\delta}\
\frac{d\tau}{\tau \Phi_p^{-1}(\tau )}\ \le\
\frac{1}{\Phi_p^{-1}(\delta_0)}\ \log\ \frac{\delta}{\delta_0}\ <\
\infty
\end{equation} because $\Phi_p^{-1}$ is non-decreasing and
$\Phi_p^{-1}(\delta_0)>0$. Moreover, by the de\-fi\-ni\-tion of the
inverse function $\Phi_p^{-1}(\tau)\equiv 0$ for all $\tau \in
[0,\tau_0],$ $\tau_0=\Phi_p(0)$, and hence (\ref{eq3.a333}) for
$\delta\in[0,\tau_0),$ generally speaking, does not imply
(\ref{eq3.333a}). If $\tau_0 > 0$, then
\begin{equation}\label{eq3.c333} \int\limits_{\delta}^{\tau_0}\
\frac{d\tau}{\tau \Phi_p^{-1}(\tau )}\ =\ \infty\ \ \ \ \ \ \ \ \ \
\forall\ \delta\ \in\ [0,\tau_0)  \end{equation} However,
(\ref{eq3.c333}) gives no information on the function $Q(z)$ itself
and, consequently, (\ref{eq3.a333}) for $\delta < \Phi(0)$ cannot
imply (\ref{eq3.333A}) at all. \end{remark}

\medskip
In view of (\ref{eq3.a333}), Theorem \ref{th5.555} follows
immediately from Lemma \ref{lem5.5C}.

\medskip
\begin{corollary} \label{cor555}{\it\, If $\Phi:[0,\infty ]\to [0,\infty ]$ is a
non-decreasing convex func\-tion and $Q$ satisfies the condition
(\ref{eq5.555}), then each of the conditions
(\ref{eq333Y})--(\ref{eq333A}) for $p\in (0, \infty)$ implies
(\ref{eq3.333A}). Moreover, if in addition $\Phi(1)<\infty$ or
$q(r)\ge 1$ on a subset of $(0,1)$ of a positive measure, then each
of the conditions (\ref{eq333Y})--(\ref{eq333A}) for $p\in (0,
\infty)$ implies
\begin{equation}\label{eq3.3} \int\limits_{0}^{1}\
\frac{dr}{rq^{\frac{\lambda}{p}}(r)}\ =\ \infty\ \ \ \ \ \ \ \ \
\forall\ \lambda\ \in\ (0,1)
\end{equation}
and also
\begin{equation}\label{eq3.3AB} \int\limits_{0}^{1}\
\frac{dr}{r^{\alpha}q^{\frac{\beta}{p}}(r)}\ =\ \infty\ \ \ \ \ \ \
\ \ \forall\ \alpha\ge 1 ,\ \beta\ \in\ (0,\alpha]\,.
\end{equation}}
\end{corollary}

\section{Sufficient conditions for equicontinuity}
\setcounter{equation}{0}

\medskip
Let $D$ be a fixed domain in the extended space
$\overline{{\C}}={\C}\cup\{\infty\}.$ Given a function $\Phi:[0,
\infty]\rightarrow [0, \infty],$ $M>0,$ $\Delta>0$,
$\frak{F}^{\Phi}_{M,\Delta}$ denotes the collection of all
homeomorphisms with finite distortion in $D$ such that
$h\left(\overline{{\C}}\setminus f(D)\right)\ge \Delta$ and
\begin{equation}\label{eq2!!}
\int\limits_D\Phi\left(K_{f}(z)\right)\frac{dx\,dy}{\left(1+|z|^2\right)^2}\
\le\ M\,.
\end{equation}

\medskip
\begin{theorem}\label{th1!}{\it\, Let
$\Phi:[0, \infty]\rightarrow [0, \infty]$ be non-decreasing convex
function. If
\begin{equation}\label{eq3!}
\int\limits_{\delta_0}^{\infty} \frac{d\tau}{\tau\Phi^{-1}(\tau)}\
=\ \infty
\end{equation}
for some $\delta_0>\tau_0:=\Phi(0),$ then the class
$\frak{F}^{\Phi}_{M,\Delta}$ is equicontinuous and, consequently,
forms a normal family of mappings for every $M\in(0, \infty)$ and
$\Delta\in(0, 1).$ }
\end{theorem}

\medskip
\begin{remark}\label{rem1}
Note that the condition
\begin{equation}\label{eq3!!}
\int\limits_D \Phi\left(K_{f}(z)\right)dx\,dy\le M
\end{equation}
implies (\ref{eq2!!}). Thus, the condition (\ref{eq2!!}) is more
general than (\ref{eq3!!}) and homeomorphisms with finite distortion
satisfying (\ref{eq3!!}) form a subclass of
$\frak{F}^{\Phi}_{M,\Delta}.$ Conversely, if the domain $D$ is
bounded, then (\ref{eq2!!}) implies the condition
\begin{equation}\label{eq4!}
\int\limits_D \Phi\left(K_{f}(z)\right)dx\,dy\le M_*
\end{equation}
where $M_*=M\cdot\left(1+\delta_*^2\right),$
$\delta_*=\sup\limits_{z\in D}|z|.$
\end{remark}

\medskip
\begin{corollary}\label{cor1!}{\,\it
Each of the conditions (\ref{eq333Y})--(\ref{eq333A}) for $p\in (0,
n-1] $ implies equicontinuity and normality of the classes
$\frak{F}^{\Phi}_{M,\Delta}$ for all $M\in (0, \infty)$ and
$\Delta\in (0, 1).$ }
\end{corollary}

\medskip
Given a function $\Phi:[0, \infty]\rightarrow [0, \infty],$ $M>0$
and $\Delta>0,$ $S^{\Phi}_{M, \Delta}$ denotes the class of all
homeomorphisms $f$ of $D$ in the Sobolev class $W_{loc}^{1,2}$ with
a locally integrable $K_{f}(z)$ such that
$h\left(\overline{{\C}}\setminus f(D)\right)\ge\Delta$ and
(\ref{eq2!!}) holds for $K_{f}(z).$ Note that if $\Phi$ is
non-decreasing, convex and non--constant on $[0,\infty)$, then
(\ref{eq2!!}) itself implies that $K_{f}(z)\in L_{loc}^1.$ Note also
that $S^{\Phi}_{M, \Delta}\subset \frak{F}^{\Phi}_{M, \Delta},$ see
e.g. Theorem 4.1 in \cite{MRSY}. Thus, we have the following
consequence.
\medskip

\begin{corollary}\label{cor2!}{\,\it
Each of the conditions (\ref{eq333Y})--(\ref{eq333A}) for $p\in (0,
1]$ implies equicontinuity and normality of the class
$S^{\Phi}_{M,\Delta}$ for all $M\in (0, \infty)$ and $\Delta\in (0,
1).$ }
\end{corollary}

\section{Necessary conditions for equicontinuity}
\setcounter{equation}{0}

\begin{theorem}\label{th3}{\it\, If the
classes $S^{\Phi}_{M,\Delta}\subset \frak{F}^{\Phi}_{M,\Delta}$ are
equicontinuous (normal) for a non--decreasing convex function
$\Phi:[0, \infty]\rightarrow [0, \infty],$ all $M\in (0,\infty)$ and
$\Delta\in (0,1).$ Then
\begin{equation}\label{eq3}
\int\limits_{\delta_*}^{\infty}\frac{d\tau}{\tau \Phi^{\,-1}(\tau)}\
=\ \infty
\end{equation}} for all $\delta_*\in (\tau_0,
\infty)$ where $\tau_0\ \colon=\ \Phi(0).$
\end{theorem}\medskip

It is evident that the function $\Phi(t)$ in Theorem \ref{th3}
cannot be constant because in the contrary case we would have no
real restrictions for $K_{I}$ except $\Phi(t)\equiv\infty$ when the
classes $S^{\Phi}_{M,\Delta}$ are empty. Moreover, by the known
criterion of convexity, see e.g. Proposition 5 in I.4.3 of
\cite{Bou}, the slope $[\Phi(t)-\Phi(0)]/t$ is nondecreasing. Hence
the proof of Theorem \ref{th3} follows from the next
statement.\medskip

\begin{lemma}\label{th3!}{\it\, Let a function $\Phi : [0,\infty]\to[0,\infty]$
be non-decreasing and %
\begin{equation}\label{eq4!!}
\Phi(t)\ \ge\ C\cdot t\qquad\forall\ t\in [T, \infty]
\end{equation}
for some $C>0$ and $T\in (0, \infty).$ If the classes
$S_{M,\Delta}^{\Phi}\subset \frak{F}_{M,\Delta}^{\Phi}$ are
equicontinuous (normal) for all $M\in (0,\infty)$ and $\Delta\in
(0,1)$, then (\ref{eq3}) holds for all $\delta_*\in (\tau_0,
\infty)$ where $\tau_0\ \colon=\ \Phi(+0).$}
\end{lemma}

\medskip
\begin{remark}\label{rem4}
Theorem \ref{th3} shows that the condition (\ref{eq3!}) in Theorem
\ref{th1!} is not only sufficient but also necessary for
equicontinuity (normality) of classes with the integral constraints
of the type either (\ref{eq2!!}) or (\ref{eq4!}) with a convex
non--decreasing $\Phi.$ In view of Proposition \ref{pr4.1aB}, the
same concerns to all the conditions (\ref{eq333Y})--(\ref{eq333A})
with $p=1.$
\end{remark}

\medskip
\begin{corollary}\label{cor3!}
{\it\, The equicontinuity (normality) of the classes
$S^{\Phi}_{M,\Delta}\subset \frak{F}^{\Phi}_{M,\Delta}$ for $M\in
(0, \infty)$, $\Delta\in (0,1)$ and non--decreasing convex $\Phi$
implies that
\begin{equation}\label{eq6!!}
\int\limits_{\delta}^{\infty}\log \Phi(t)\ \frac{dt}{t^{2}}\ =\
\infty
\end{equation}
for all $\delta>t_0$ where $t_0:=\sup\limits_{\Phi(t)=0}t,$ $t_0=0$
if $\Phi(0)>0.$
 }
\end{corollary}

\medskip
The condition (\ref{eq6!!}) is also sufficient for
equi\-con\-ti\-nui\-ty (normality) of the classes
$S^{\Phi}_{M,\Delta}$ and $\frak{F}^{\Phi}_{M,\Delta}$.

\end{document}